\documentclass[11pt]{amsart}
\usepackage{amsmath} \usepackage{amssymb}
\theoremstyle{plain} \newtheorem{theor}{Theorem}
{}
\newtheorem{lemma}{Lemma}

\newtheorem{corol}{Corollary} {}
\theoremstyle{definition} \newtheorem{defin}{Definition}

\theoremstyle{remark}  
\numberwithin{equation}{section}
\DeclareMathOperator{\cov}{cov} 
 
 \DeclareMathOperator{\cof}{cof}

\newcommand{\IF}{\text{ if }} 
\newcommand{\AND}{\text{ and }}  \newcommand{\forces}[2]{\Vdash_{#1} \mbox{``} #2 \mbox{''}}

\newcommand{\Reals}{{\mathbb R}}
\newcommand{\Rationals}{{\mathbb Q}}
\newcommand{\Irrationals}{\Reals\setminus{\mathbb Q}}
\newcommand{\Poset}{{\mathbb P}}

\newcommand{\Integers}{{\mathbb Z}}

\newcommand{\presup}[2]{\, ^{#1} \! #2}
\newcommand{\wpresup}[1]{\presup{\stackrel{\omega}{\smile}}{#1}}
\newcommand{\fomom}{\presup{\omega}{\omega}}
\newcommand{\wfomom}{\wpresup{\omega}}
  \newcommand{\card}[1]{\lvert #1
\rvert}
\newcommand{\cccc}[2]{#1{}^\wedge #2}

\title[Covering Numbers of Trees]{Covering Numbers Associated with
Trees Branching into a Countably Generated Set of Possibilities} \author[S.
Shelah]{Saharon Shelah} \address{
Institute of Mathematics, Hebrew University}
\date{}
\thanks{The author would like to thank Juris Stepr\={a}ns for writing
up this paper. This paper is number 660 in the author's personal listing.}
\begin{document}
%
\maketitle \tableofcontents \bibliographystyle{plain}
\section{Introduction}
This paper is concerned with certain generalizations of meagreness and
their combinatorial equivalents. The simplest example, and the one
which motivated further study in this area, comes about by
considering the following definition:
\begin{defin} For any $A\subseteq \Reals$ a set $X\subseteq \Reals$
will be said to be  $A$-{\em  nowhere
dense}\footnote{This notation corrects the terminology of
\cite{step.30} which called a special case of this notion {\em almost nowhere
dense} in spite of the fact that almost nowhere dense sets are nowhere
dense rather than the converse.} if and
only if for every  $q \in A$ there exists and integer $k \in
\Integers$ such that the interval whose endpoints are $q$ and $q+ 1/k$
is disjoint from $X$. A set which is the union of countably many
$A$-nowhere dense sets will be called $A$-{\em very meagre}.
\end{defin}

The notion of an $A$-nowhere dense set for various subsets $A$ of the
reals may prove to be of interest in its own right, but this paper will
be concerned exclusively with the special case $A=\Rationals$.
Notice that rational perfect
sets introduced by Miller in \cite{mill.ratpfct} form a subset of the
$\Rationals$-nowhere dense sets since the closure of a set is rational
perfect if it is
perfect and disjoint from the rationals. On the other hand, a set is
$\Rationals$-nowhere dense if its  2-sided-closure is disjoint
from the rationals where the 2-sided-closure of a set refers to all
those reals which are limits of both decreasing and increasing
sequences from the set.

In \cite{step.30} the least
number of $\Rationals$-meagre sets required to cover the real line is
examined and is denoted by ${\mathfrak d}_1$. It is shown in
\cite{step.38} that there is
a continuous function $H$ --- first constructed by
Lebesgue --- such that  the least number of smooth functions into which
$H$ can be decomposed is equal to ${\mathfrak d}_1$.
This paper will further study ${\mathfrak d}_1$ and some of its
generalizations. As well, an  equivalence will be established between
 $\Rationals$-meagreness and certain combinatorial
properties of trees. This will lead to new cardinal invariants and
various independence results about these  will then be established.

\section{Equivalences}
 \begin{defin}
A set $X\subseteq \presup{\omega}{\Integers}$ is small if for each
$n\in\omega$ and $\sigma: n \to \Integers$ there is  some $k\in
\Integers$ such that either $$\{f(n) : f \in X \AND \sigma \subseteq
f\} \subseteq \{m\in \Integers : m < k\}$$
or
 $$\{f(n) : f \in X \AND \sigma \subseteq
f\} \subseteq \{m\in \Integers : m > k\}$$
\end{defin}
It will be shown that the least number of very meagre sets required to
cover $\Reals$ is equal to least number of small sets required to cover
$\presup{\omega}{\Integers}$. This is an immediate consequence of the
following lemma. The idea is to modify one of the
standard\footnote{For example see page 5 of \cite{miller.dst}.}
constructions of a homeomorphism between $\fomom$ and $\Reals\setminus
\Rationals$. In this construction
$\fomom$ is identified with $\presup{\omega}{{\mathbb Z}}$ and finite
sequences $\sigma : m \to {\mathbb Z}$ are mapped to open intervals $I(\sigma)$
so that
\begin{itemize}
\item if $\sigma\subseteq \tau$ then $I(\sigma)\supseteq I(\tau)$
\item $\bigcup_{n\in {\mathbb Z}}\overline{I(\cccc{\sigma}{n})} =
I(\sigma)$
\item the left endpoint of $I(\cccc{\sigma}{n})$ is the right endpoint
of $I(\cccc{\sigma }{(n-1)})$.
\end{itemize}
This construction will be modified by, essentially, mapping
sequences which end with a negative integer to their counterpart above the
neighbour to the right of their parent sequence. The details are
explained in the following.
\begin{lemma}
\label{l:wu}
Let $X\subseteq\presup{\omega}{\Integers}$ be the set of all sequences
eventually equal to $-1$ together with the constant sequence $0$.  There is a
bijection $F:\presup{\omega}{\Integers}\setminus X\to \Irrationals$ such that
$A\subseteq \Irrationals$ is very nowhere dense if and only if $F^{-1}A$ is
small.
\end{lemma}
\begin{proof} To begin, let $\{q_n\}_{n\in\omega}$ enumerate $\Rationals$ so
that each rational occurs infinitely often.  Let $\wpresup{\Integers}$ denote
the set of finite sequences of integers --- in other words, the set of
functions from an integer to $\Integers$. For $\sigma\in\wpresup{\Integers}$
define $\mu(\sigma) = \sigma(n)$ where $n $ is the greatest integer in the
domain of $\sigma$.  Let $\prec$ be the partial ordering of
$\wpresup{\Integers}$ defined by $\sigma \prec \tau$ if 
\begin{itemize}
\item there exists some least integer $n$ such that
$\sigma(n)  \neq \tau(n)$ and $0 \leq \tau(n) < \sigma(n)$
\item there exists some least integer $n$ such that
$\sigma(n)  \neq \tau(n)$ and $\tau(n) < \sigma(n) < 0$
\item there exists some least integer $n$ such that
$\sigma(n)  \neq \tau(n)$ and $ \sigma(n) < 0 \leq \tau(n)$
\end{itemize}
 For $\sigma \in
\wpresup{\Integers}$ and $n\in \Integers$ define $\sigma^n \in
\wpresup{\Integers}$
such that $\sigma$ and $\sigma^n$ have the same domain, $\sigma(i) =
\sigma^n(i)$ if $i$ is not the maximal element of the common domain of
$\sigma$ and $\sigma^n$ and $\mu(\sigma)+n=\mu(\sigma^n)$. Moreover, let
$\sigma^-=\sigma\restriction(\card{\sigma}-1)$. Next, construct by induction
a mapping $f:\wpresup{\Integers} \to \Reals$ such that:
\begin{itemize}
\item $f$ is a $\prec$ order preserving mapping whose range is
disjoint from $\Rationals$
\item if $\mu(\sigma)\neq -1$ and $\sigma^-$ is $-1$ on its domain, or
$\sigma^-\not\equiv 0$, $\sigma^-\not\equiv -1$, $\mu(\sigma)\geq 0$, then
$f(\sigma^{1})< f(\cccc{\sigma}{i})<f(\sigma)$ for each $i\in\Integers
\setminus\{0,-1\}$, and $f(\cccc{\sigma}{-1})=f(\sigma^1)$, $f(
\cccc{\sigma}{0})=f(\sigma)$
\item if $\mu(\sigma)\neq0$ and $\sigma^-$ is $0$ on its domain, or $\sigma^-
\not\equiv 0$, $\sigma^-\not\equiv-1$, $\mu(\sigma)<0$, then $f(\sigma)<
f(\cccc{\sigma}{i})<f(\sigma^{-1})$ for each $i\in\Integers\setminus\{0,
-1\}$, and $f(\cccc{\sigma}{-1})=f(\sigma)$, $f(\cccc{\sigma}{0})=f(
\sigma^{-1})$
\item $\lim_{n\in\omega}f(\sigma^{n}) = \lim_{n\in\omega}f(\sigma^{-n})\in
\Rationals$ for every $\sigma\in \wpresup{\Integers}$
\item $\card{f(\sigma)-f(\tau)}<\frac{1}{k + 1}$ for every $\sigma$ and $\tau$
in $\presup{k}{\Integers}$ such that $\tau$ is the immediate successor of
$\sigma$ with respect to $\prec$ in $\presup{k}{\Integers}$
\item if $\sigma\equiv 0$ then $f(\cccc{\sigma}{0})=f(\sigma)+\frac{1}{
\card{\sigma}+1}$, $f(\cccc{\sigma}{-1})=f(\sigma)$
\item if $\sigma\equiv -1$ then $f(\cccc{\sigma}{-1})=f(\sigma)-\frac{1}{
\card{\sigma}+1}$, $f(\cccc{\sigma}{0})=f(\sigma)$
\item if $\sigma,\tau\in\presup{k}{\Integers}$ are two successive sequences
(with respect to $\prec$), and both $q_k$ and $f(\cccc{\sigma}{1})$ are
between $f(\sigma)$ and $f(\tau)$ then $\lim_{n\in\omega}f(\cccc{\sigma}{
{n}})=\lim_{n\in\omega}f(\cccc{\sigma}{-n})=q_k$
\end{itemize}
For $\sigma\in\wpresup{\Integers}$ define an interval of reals by $I(\sigma)
=[f(\sigma),f(\sigma^{*})]$, where $\sigma^*$ is suitably either successor or
predecessor of $\sigma$ in $\presup{\card{\sigma}\;}{\Integers}$, and define
$F:\presup{\omega}{\Integers}\setminus X\to \Irrationals$ by taking the
intersection along a branch --- in other words, $F(\sigma) $ is the unique
element of $\bigcap_{n\in\omega}I(\sigma\restriction n)$.
It is easy to check that this mapping has the desired properties.
\end{proof}

\section{Trees of Countable Structures}
A cover on a countable set $X$ is a countable subset ${\mathcal B}
\subseteq {\mathcal P}(X)$ such that
\begin{itemize}
\item $X \notin \mathcal B$
\item if $B \in \mathcal B$ and $b \in [X]^{<\aleph_0}$
then $B\cup b\in\mathcal B$.
\end{itemize}
 For a
cover $\mathcal B$ on $X$ define $\overline{\mathcal B} = \{Y\subseteq
X : (\exists A\in {\mathcal B})(Y\subseteq A)\}$ and define ${\mathcal
B}^+ = {\mathcal P}(X) \setminus \overline{\mathcal B}$. Define a
$\mathcal B$-tree to be a tree $T\subseteq
\presup{\stackrel{\omega}{\smile}}{X}$ such that for each $t\in T$ the
set of successors of $t$ in $T$ belongs to $\overline{\mathcal B}$ or,
to be more precise, $\{s(\card{t}) : s\in T \AND t\subsetneq s\}
\in \overline{\mathcal B}$.
The notation $\wfomom$ will be used to denote the set of all functions
from a proper, initial segment of $\omega$ to $\omega$.
Finally, define ${\mathcal J}_{\mathcal
B}$ to be the ideal generated by all sets $X\subseteq \wfomom$ such that
there is
some ${\mathcal B}$-tree $T$ such that $X \subseteq \overline{T} = \{f
: (\forall n\in\omega)f\restriction n \in T\}$.  Note that ${\mathcal
J}_{\mathcal B}$ is a countably complete ideal.

The examples of covers with which this paper will be concerned are of
the form
$${\mathcal B}_n = \{A\subseteq \omega\times n : (\exists i \in
n)\card{\{m\in\omega : (m,i) \in A\}} < \aleph_0\}$$ although many
other examples are possible. It will be shown that for any integer $m
> 1$
 it is consistent that $\cov({\mathcal
J}_{{\mathcal B}_{m}}) = \omega_2$ but
$\cov({\mathcal
J}_{{\mathcal B}_{m+1}}) = \omega_1$.

\begin{defin}\label{d:poset}
For a cover ${\mathcal B}$ on a set $X$ define $\Poset({\mathcal B})$
to be the set of all triples $(t, F, \Gamma)$ such that:
\begin{enumerate}
\item
$t \subseteq \presup{\stackrel{\omega}{\smile}}{X}$ is a finite
subtree --- in particular, $t$ is closed under initial segments
\item
$F:t\to {\mathcal B}$
\item
$\Gamma \in [\presup{\omega}{X}]^{<\aleph_0}$
\item there is a one-to-one function $\beta : \Gamma \to t$ which maps
$\Gamma$ onto the maximal
nodes of $t$ such that $\beta(x)\subseteq x$ for all $ x\in \Gamma$
\item
$x(n) \in F(x\restriction n)$ for every $x\in \Gamma$ and for every
$n\in \omega$ such that $x\restriction n \in t$
\end{enumerate}
The ordering on $\Poset({\mathcal B})$ is coordinatewise
containment. Observe that if $(t,F,\Gamma)$ and
$(t',F',\Gamma')$ are in ${\mathbb P}({\mathcal B})$ and $t=t'$, $F=F'$
then the two
conditions are compatible.
\end{defin}

If ${\mathcal B}$ is a cover on $X$ and ${\mathcal C}$ is a cover on $Y$
then define ${\mathcal B} \prec {\mathcal C}$ if and only if for every
$A \in {\mathcal B}^+$ and  for every $H:A \to \presup{\omega}{Y}$
there is $B\subseteq A$ such that $B \in
{\mathcal B}^+$ such that
 there is a finite  $t\subseteq \presup{\stackrel{\omega}{\smile}}{Y}$
and a mapping $F:t\to {\mathcal C}$ and
  there is a finite set $C \subseteq \presup{\omega}{Y}$ and $B_0\cup
B_1 \subseteq B$ such that $B_0\cup B_1 \in {\mathcal B}^+$ and:
\begin{equation}\label{e:split}
(\forall \{b, b'\} \in [B_0]^2)(\exists \tau \in t)(\exists \{y, y'\}
\in [F(\tau)]^2)(\cccc{\tau}{y} \subseteq H(b) \AND
\cccc{\tau}{y'} \subseteq H(b'))\end{equation}
\begin{equation}\label{e:seq}
(\forall n\in \omega)(\card{\{b\in B_1 : (\forall c \in C)( H(b)\restriction
n \not\subseteq c)\}}< \aleph_0)\end{equation}

\begin{defin} \label{Lusin}
A $({\mathcal J},\kappa)$-Lusin set for an ideal $\mathcal J$ on a set $X$
is a set
$L\subseteq X$ such that $\card{L\cap J} < \kappa$ for all $J\in
\mathcal J$.  \end{defin}

\begin{lemma} If $\mathcal B$ \label{l:m}and $\mathcal C$ are covers on $X$ and
$Y$ respectively, ${\mathcal B } \prec \mathcal C$, $\kappa$ is a
cardinal of uncountable cofinality and $W$ is a $({{\mathcal J}_{\mathcal
B}}, \kappa)$-Lusin set
then
$1\forces{\Poset({\mathcal C})}{
W \text{ is a } ({{\mathcal J}_{\mathcal B}}, \kappa)\text{-Lusin
set}}$.\end{lemma}
\begin{proof}
If not, then let $T$ be a
$\Poset({\mathcal C})$-name for ${\mathcal B}$-tree such that
$$1\forces{\Poset({\mathcal C})}{\{w_\alpha\}_{\alpha \in
\kappa}\subseteq \overline{T}\cap W \AND w_\alpha\neq w_\beta \IF
\alpha\neq\beta}$$ and,
for each $\alpha \in \kappa$, choose $(t_\alpha, F_\alpha,
\Gamma_\alpha) \in \Poset({\mathcal C})$ deciding the value of
$w_\alpha$ --- in other words,  such that there is ${y}_\alpha$ such
that
 $(t_\alpha,
F_\alpha, \Gamma_\alpha) \forces{ \Poset({\mathcal C})}{w_\alpha =
\check{y}_\alpha}$. Then, using the fact that $\kappa$ has uncountable
cofinality,
 choose $t$ and $F$ such that $W^* = \{ y_\alpha :
 t= t_\alpha \AND F =  F_\alpha\}$ has cardinality $\kappa$. It
follows $W^*$ is also a $({{\mathcal J}_{\mathcal B}}, \kappa)$-Lusin
set and, hence,  that it is in ${\mathcal J}_{\mathcal B}^+$.

Next, choose $m\in\omega$ and $\sigma : m\to\omega$ such that
$$S_0 = \{s\in X : (\exists w \in W^*)(\cccc{\sigma}{s}\subseteq w)\} \in
 {\mathcal B}^+$$ and choose $W'\subseteq W^*$ such that for each $s
\in S_0$ there is a unique $w_s \in W'$ such that $w_s(m) = s$. Let
$\alpha(s)$ be the unique ordinal such that $w_s = y_{\alpha(s)}$.
Let the maximal nodes of $t$ be enumerated by $\{\nu_i\}_{i=1}^ k$ and
let $b^s_i$ be the unique member of $\Gamma_{\alpha(s)}$ such that
$\nu_i\subseteq b^s_i$.
Proceed by induction to define $S_0^i$ and $S_1^i$ for $i\leq k$ so
that if
$S_0^i$ and $S_1^i$ have been chosen then $S_0^{i+1}\cup
S_1^{i+1}\subseteq S_i=S_0^i \cup S_1^i$ is
chosen so that
 $S_0^{i+1}\cup S_1^{i+1} \in {\mathcal B}^+$
and so that
there is a finite  $r_{i+1}\subseteq \presup{\stackrel{\omega}{\smile}}{Y}$
and a mapping $F_{i+1}:r_{i+1}\to {\mathcal C}$ and
  there is a finite set $L_{i+1} \subseteq \presup{\omega}{Y}$ such that
\begin{equation}\label{e:spliti}
(\forall \{s,\bar{s} \} \in [S_0^{i+1}]^2)(\exists \tau \in
r_{i+1})(\exists \{y, \bar{y}\}
\in [F_{i+1}(\tau)]^2)(\cccc{\tau}{y} \subseteq b^{s}_{i+1} \AND
\cccc{\tau}{\bar{y}} \subseteq {b}^{{\bar{s}}}_{i+1})\end{equation}
\begin{equation}\label{e:seqi}
(\forall n\in \omega)(\card{\{s\in S_1^{i+1} : (\forall c \in L_{i+1})
(b^{s}_{i+1}\restriction
n \not\subseteq c)\}}< \aleph_0)\end{equation}
\begin{equation}\label{e:a}
\nu_i \subseteq \bigcap( r_i\cup L_i)
\end{equation}
\begin{equation}\label{e:b}
F_i(\nu_i)=F(\nu_i)
\end{equation}
\begin{equation}
(\forall c \in L_i)(\exists \sigma \in \max(r_i))(\sigma \subseteq c)
\end{equation}
\begin{equation}
(\forall c \in L_i)(\forall \sigma \in r_i)(\sigma \subseteq c
\rightarrow c(|\sigma|) \in F(\sigma))
\end{equation}
\begin{equation}
\IF \tau\in \Gamma_{i+1} \AND n<|\tau| \text{ then } \tau(n)\in
F_{i+1}(\tau\restriction n)\end{equation}
This is easily accomplished using the definition of ${\mathcal B}
\prec {\mathcal C}$. Then define a new condition $(t', F', \Gamma')$
such that
$t\cup \bigcup_{i=1}^kr_i\subseteq t'$, $F\cup
\bigcup_{i=1}^kF_i\subseteq F'$
and $\Gamma\cup \bigcup_{i=1}^kL_i\subseteq
\Gamma'$.

It suffices to show that $(t', F', \Gamma')\forces{\Poset({\mathcal
C})}{\{s\in S_k : w_s \in \overline{T}\}\notin
\overline{\mathcal B}}$ because this would contradict that
$1\forces{\Poset({\mathcal
C})}{T \text{ is a ${\mathcal B}$-tree}}$. Therefore suppose that
$(t'', F'', \Gamma'')$ is a condition
extending $ (t', F', \Gamma')$ such that $$(t'', F'',
\Gamma'')\forces{\Poset({\mathcal
C})}{ {\{s\in S_k : w_s \in \overline{T}\}\subseteq  B}}$$ for
some $B \in \mathcal B$.  In order to obtain a contradiction, notice
that for each $i$ the set  $\Lambda_0^i = \{s
\in S_k\cap S_0^i : (\exists \tau \in t'')(\tau
\subseteq b^{s}_i \AND  b^{s}_i(\card{\tau})\notin F_i''(\tau))\}$ is
finite because
of \eqref{e:spliti} and the fact that
$r_i \subseteq t''$. On the
other hand, it is possible to choose $N$ so large that for each $i$
and each  $c \in L_i$  the sequence  $c\restriction N$ does not
belong to $t''$. Then for each $i$ the set
$$\Lambda_1^i = \{s \in S_k\cap S_1^i : (\forall c\in
L_i)(b^{s}_i\restriction N
\not\subseteq c)\}$$ is also finite.
Therefore it is possible to
choose $$ s\in S_k \setminus (B \cup \bigcup_{i=1}^k(\Lambda_0^i \cup
\Lambda_1^i) )$$
and let $\Gamma^* = \Gamma'' \cup \{b^{s}_i\}_{i=1}^k$.

First note
that it is easy to
extend $t'' $ to $t^*$ and $F''$ to $F^*$ so that $(t^*,
F^*, \Gamma^*)$ satisfies Conditions~1, 2 and 4 of
Definition~\ref{d:poset} because $[Y]^{<\aleph_0} \subseteq \mathcal
C$.
 Since Condition~3 is also satisfied
it  suffices
to show that $(t^*,
F^*, \Gamma^*)$ satisfies Condition~5 in Definition~\ref{d:poset}.
 This, in turn, follows
from consideration of the two cases. First,
if $s \in S_0^i$ then $s \notin \Lambda_0^i$ and so
there is some $\tau \in r_i$ and $ y \in
F_i(\tau) = F''(\tau) = F^*(\tau)$ such that $\cccc{\tau}{y}
\subseteq b^{s}_i$ and $ \cccc{\tau}{y} \notin t''$. Hence there
is no contradiction to Condition~5 because any such contradiction
would already have occured in the condition $(t_{\alpha(s)}, F_{\alpha(s)},
\Gamma_{\alpha(s)})$.
If $s \in S_1^i$ then there is some $c \in L_i\subseteq \Gamma''$ such that
$b^{s}_i\restriction N = c\restriction N$ and $c$ does not violate
Condition~5. Moreover the choice of $N$ guarantees that $c$  violates
Condition~5 if and only if $b^{s}_i$ does.
But now $(t^*,
F^*, \Gamma^*)$ obviously extends $(t,F,\Gamma_{\alpha(s)}) =
(t_{\alpha(s)},F_{\alpha(s)},\Gamma_{\alpha(s)}) $ which  forces that
$w_s \in \overline{T}$ and, hence, that $s\in B$ and this
is a contradiction.
\end{proof}

\begin{corol}
Suppose that  \label{c:1}$\mathcal B$ is a cover and that  $W$ is a
$({{\mathcal J}_{\mathcal B}}, \kappa)$-Lusin set. Suppose also that
$\lambda \in \kappa$ and $\cof(\kappa) > \aleph_0$
 and that  $\{{\mathcal C}_\alpha\}_{\alpha \in
\lambda}$ are covers such that ${\mathcal B} \prec {\mathcal
C}_{\alpha}$ for each $\alpha \in \lambda$. If, furthermore,  $\Poset$ is the
finite support iteration of $\{\Poset({\mathcal C}_\alpha)\}_{\alpha \in
\lambda}$ then $1\forces{\Poset}{ W \text{ is a }
({{\mathcal J}_{\mathcal B}}, \kappa)\text{-Lusin set}}$.\end{corol}
\begin{proof} Proceed by induction on $\lambda$. If $\lambda$ is a limit
ordinal and the lemma fails  then let $G\subset \Poset$ be generic
and choose a ${{\mathcal J}_{\mathcal B}}$-tree $T$ and
$\{p_\alpha\}_{\alpha\in\kappa}\subseteq G$ such that
$p_\alpha\forces{\Poset}{\hat{w}_\alpha \in \overline{T}\cap W}$ and such
that the $w_\alpha$ are all distinct.
Then use the finite support of the
iteration and the fact that
the cofinality of $\lambda$ is less than $\kappa$ to conclude that
there is some $\beta \in \alpha$ such that $G\cap
\Poset_\beta$ --- where $\Poset_\beta$ is the
finite support iteration of $\{\Poset({\mathcal C}_\nu)\}_{\nu \in
\beta}$ --- contains $\kappa$ of the conditions
$\{p_\alpha\}_{\alpha\in\kappa}$. This contradicts the
induction hypothesis because the closure of the corresponding
$w_\alpha$ will form a ${\mathcal J}_{\mathcal B}$-tree contained in
$\overline{T}$.

At successors, use the induction hypothesis,
Lemma~\ref{l:m} and the absoluteness of the relation $\prec$.\end{proof}

\begin{lemma}
If \label{l:st}$\Poset$ is the finite support iteration of length $\omega$
of the
partial orders $\Poset({\mathcal B})$ for some cover $\mathcal B$ and
$G$ is $\Poset$ generic over $V$ then, in $V[G]$, there are countably
many sets in ${\mathcal J}_
{\mathcal B}$ whose union covers $V\cap
\fomom$.
\end{lemma}
\begin{proof} Standard. \end{proof}

\begin{theor}
If $\kappa_0 > \kappa_1$ are uncountable regular
cardinals and ${\mathcal B}_0$ and ${\mathcal B}_1$ are covers such
that ${\mathcal B}_0 \prec {\mathcal B}_1$ \label{t:m}
then it is consistent that $\cov({\mathcal J}_{{\mathcal B}_0}) =
\kappa_0$ and  $\cov({\mathcal J}_{{\mathcal B}_1}) =
\kappa_1$
\end{theor}
\begin{proof} Let $V\models$ GCH \& ZFC and let $\Poset_0$ be Cohen forcing
for adding $\kappa_0$ Cohen reals and let $\Poset_1$ be
the finite
support iteration, of length $\kappa_1$ of the partial orders
$\Poset({\mathcal B}_1)$. Let $\Poset = \Poset_0 * \Poset_1 $
and let $G$ be $\Poset$ generic over $V$

Notice first that,
since it is easily verified that each member of any
${\mathcal J}_{{\mathcal B}}$ is meagre, it follows that the
$\kappa_0$-Lusin set added by $\Poset_0$
is also an
$({{\mathcal J}_{{\mathcal B}_0}}, \kappa_0)$-Lusin set.
\relax From Corollary~\ref{c:1} it follows that $V[G]$ has a
$({\mathcal J}_{{\mathcal B}_0}, \kappa_0)$-Lusin set.
Hence $\cov({\mathcal J}_{{\mathcal B}_0})\geq \kappa_0$. Since the
continuum in $V[G]$ is $\kappa_0$ it follows that $\cov({\mathcal
J}_{{\mathcal B}_0}) = \kappa_0$ in this model.
On the other hand, it follows from Lemma~\ref{l:st} that
$\cov({\mathcal J}_{{\mathcal B}_1}) \leq \kappa_1$. Observe that the
finite support iteration $\Poset_1$ adds Cohen reals over each
intermediate model. Moreover, as it has already been observed that
 each member of any
${\mathcal J}_{{\mathcal B}}$ is meagre, it follows no family of size
less than $\kappa_1$ elements of ${\mathcal J}_{{\mathcal B}_1}$
can cover all the reals. Hence $\cov({\mathcal
J}_{{\mathcal B}_1}) = \kappa_1$ in this model.
\end{proof}

\begin{defin}\label{d:ex}
Define ${\mathcal B}^i_n$ to consist of all sets $$B(a,k) =
(\omega \times a) \cup (k\times  (n\setminus a))$$ where $ k \in \omega$
and $a \in [n]^i$.\end{defin} It is easy to check that each ${\mathcal
B}^i_n$ is a cover on $\omega \times n$ if $ i < n$.

\begin{lemma}\label{l:pr} If  $i + 1 \leq j$ then
${\mathcal B}_n^i \prec {\mathcal B}_m^{j}$.
\end{lemma}
\begin{proof}
Let $A \in ({\mathcal B}^i_n)^+$ and suppose that
$H:A \to \presup{\omega}{(\omega\times m)}$.
Without loss of generality it may be assumed that $A = \bigcup_{k\in
a}A_k\times \{k\}$ where $a \in [n]^{i+1}$ and each $A_k$ is
infinite. Let $D$  be the set of all
$k\in a $ such that the range of $H\restriction A_k\times\{k\}$ is
contained in a finite branching tree. By compactness, for each $k \in
D$ there is $c_k \in \presup{\omega}{(\omega\times m)}$ and $B^k\in
[A_k\times\{k\}]^{\aleph_0}$ such that for
all $N\in \omega$ the set $\{x\in B^k : H(x)\restriction N
\not\subseteq c_k\}$ is finite.

For $k \in a \setminus D$ there must exist some $t_k \in
\presup{\stackrel{\omega}{\smile}}{(\omega\times m)}$ such that
$$\{z\in
\omega\times m : (\exists x\in A_k\times\{k\})( \cccc{t_k}{z}\subseteq
H(x))\}$$ is infinite. It is then possible to choose $J(k) \in m$ and
$B^k\in [A_k\times\{k\}]^{\aleph_0}$ such that $$\{z\in
\omega\times m : (\exists x\in B^k)( \cccc{t_k}{z}\subseteq
H(x))\}\subseteq \omega\times\{J(k)\}$$ and, if $\{x,y\} \in [B^k]^2$
and $H(x) = \cccc{t_k}{x'} $ and $H(y) = \cccc{t_k}{y'}$ then $x' \neq
y'$. Now let $B_0 = \bigcup_{k\in a\setminus D}B^k$ and $t =
\{t_k : k\in a\setminus D\}$. For $s\in t$ let
 $$F(s) = \cup\{\omega\times\{J(k)\} : t_k = s\}$$ and  let $B_1 =
\bigcup_{k\in  D}B^k$ and
 $C =
\{c_k : k\in a\setminus D\}$.
After shrinking $B^k$ further to guarantee that
if $k_0,k_1\in a\setminus D$, $t_{k_0}=t_{k_1}$,
  $J(k_0)=J(k_1)$ and $x\in B^{k_0}$, $y\in B^{k_1}$
then $H(x)(|t_{k_0}|)\neq H(y)(|t_{k_0}|)$,
it is routine to check that these
choices  witness the required instance of ${\mathcal B}_n^i \prec
{\mathcal B}_m^{j}$.
\end{proof}

\begin{lemma}\label{l:9}
If $j < n \leq m$ then $\cov({\mathcal J}_{{\mathcal B}_n^j}) =
\cov({\mathcal J}_{{\mathcal B}_m^j})$.
\end{lemma}
\begin{proof}
Let $\beta_{m,n} : \omega\times m \to \omega\times n$ be  a bijection which is
the identity on $\omega\times (n-1)$. This induces a bijection
$$\bar{\beta}_{m,n} : \presup{\omega}{(\omega\times m)}\to
\presup{\omega}{(\omega\times n)}$$
which sends members of the ideal ${\mathcal J}_{{\mathcal B}_m^j}$ to
members of the ideal ${\mathcal J}_{{\mathcal B}_n^j}$. Hence
$\cov({\mathcal J}_{{\mathcal B}_n^j}) \leq
\cov({\mathcal J}_{{\mathcal B}_m^j})$.

To prove the other inequality, proceed by induction on $n\geq j+1$ to
show that $\cov({\mathcal J}_{{\mathcal B}_n^j}) \geq
\cov({\mathcal J}_{{\mathcal B}_{n+1}^j})$.
Now, let ${\mathcal T}$ be a family of no more than $\cov({\mathcal
J}_{{\mathcal B}_n^j})$ many ${\mathcal B}_n^j$-trees such that
$\presup{\omega}{(\omega\times n)} = \bigcup_{T\in {\mathcal
T}}\overline{T}$.
It follows that $\bar{\beta}_{n+1,n}^{-1}(\overline{T})$ is the
closure of a ${\mathcal B}_{n+1}^{j+1}$-tree. Next,
 note that if $T$ is a ${\mathcal B}_{n+1}^{j+1}$-tree then, by
identifying the successor nodes of $T$ with the elements of
$\omega\times (j+1)$, it is possible to construct a natural
bijection $\gamma_T$ from $T$ to
$\presup{\stackrel{\omega}{\smile}}{(\omega\times
(j+1))}$ which also induces a bijection $\bar{\gamma}_T : \overline{T}
\to\presup{\omega}{(\omega\times (j+1}))$.
Furthermore, if $A\in
{\mathcal J}_{{\mathcal B}_{j+1}^j}$ then
$\bar{\gamma}_T^{-1}(A) \in {\mathcal J}_{{\mathcal B}_{n+1}^j}$.
  By the induction hypothesis it is known that $\cov({\mathcal
J}_{{\mathcal B}_n^j}) = \cov({\mathcal
J}_{{\mathcal B}_{j+1}^j})$ and so it is possible to find
a family ${\mathcal S}_T$ of no more than $\cov({\mathcal
J}_{{\mathcal B}_n^j})$ many ${\mathcal B}_{j+1}^j$-trees such that
$\presup{\omega}{(\omega\times (j+1))} = \bigcup_{T\in {\mathcal
S}}\overline{T}$. Now, for $T\in \mathcal T$ and $S \in \mathcal S$
let $B(T,S) = \bar{\gamma}_T^{-1}\overline{S}$ and note that $B(T,S)$
is contained in the closure of a ${\mathcal B}_{n+1}^{j}$-tree. Hence
$\{B(T,S) : T\in \mathcal T \AND S\in {\mathcal S}\}$ witnesses that
$\cov({\mathcal J}_{{\mathcal B}_n^j}) \geq
\cov ({\mathcal J}_{{\mathcal B}_{n+1}^j})$.
\end{proof}

\begin{corol}If $j < n $, $i < m$
and  $j > i$.
then it is consistent that
$\cov({\mathcal J}_{{\mathcal B}_n^j}) <  \cov({\mathcal J}_{{\mathcal
B}_m^i})$.
\end{corol}
\begin{proof}
Apply Lemma~\ref{l:pr} and Lemma~\ref{l:9} together with
Theorem~\ref{t:m}.
\end{proof}

As a final remark,
 notice that  an immediate consequence  of Theorem~\ref{t:m} and
 Lemma~\ref{l:wu}  is
that, for any pair of uncountable, regular cardinals $\kappa_0 >
 \kappa_1$,
 it is consistent that $\cov({\mathcal M}) = \kappa_0$ and
${\mathfrak d}_1 = \kappa_1$ or, in other words, the covering number
 of the meagre ideal is $\kappa_0$ while the covering number of the
 $\Rationals$-meagre ideal is $\kappa_1$.


\begin{thebibliography}{1}

\bibitem{mill.ratpfct}
A.~Miller.
\newblock Rational perfect set forcing.
\newblock In D.~A.~Martin J.~Baumgartner and S.~Shelah, editors, {\em Axiomatic
  Set Theory}, volume~31 of {\em Contemporary Mathematics}, pages 143--159,
  Providence, 1984. American Mathematical Society.

\bibitem{miller.dst}
Arnold~W. Miller.
\newblock {\em Descriptive Set Theory and Forcing}, volume~4 of {\em Lecture
  Note in Logic}.
\newblock Springer, Berlin, 1995.

\bibitem{step.38}
J.~Stepr\={a}ns.
\newblock Decomposing with smooth sets.
\newblock {\em Trans. Amer. Math. Soc.}, 0:0--0, 1900.

\bibitem{step.30}
J.~Stepr\={a}ns.
\newblock A very discontinuous {B}orel function.
\newblock {\em J. Symbolic Logic}, 58:1268--1283, 1993.

\end{thebibliography}
\end{document}